\documentclass{elsart}
\usepackage{latexsym}
\usepackage{amssymb}
\usepackage{amsmath}
\usepackage{amsfonts}
\usepackage{amsopn}
\usepackage{amscd}
\usepackage{graphicx}

\newtheorem{theorem}{Theorem}
\newtheorem{lemma}{Lemma}

\begin{document}

\begin{frontmatter}
\title{Infinite permutations of lowest maximal pattern complexity}
\author[N]{S. V. Avgustinovich}
\ead{avgust@math.nsc.ru}
\author[N]{A. Frid}
\ead{anna.e.frid@gmail.com}
\author[M]{T. Kamae}
\ead{kamae@apost.plala.or.jp}
\author[N]{P. Salimov}
\ead{ch.cat.s.smile@gmail.com}
\address[N]{Sobolev Institute of Mathematics SB RAS
\\ Koptyug av., 4, 630090 Novosibirsk, Russia} 
\address[M]{Matsuyama University, 790-8578 Japan}
\begin{abstract}
An infinite permutation $\alpha$ is a linear ordering of $\mathbb N$.
We study properties of infinite permutations analogous to those of infinite words,
and show some resemblances and some differences between permutations and words.
In this paper, we define maximal pattern complexity $p^*_{\alpha}(n)$ for
infinite permutations and show that this complexity function is ultimately
constant if and only if the permutation is ultimately periodic; otherwise
its maximal pattern complexity is at least $n$, and the value
$p^*_{\alpha}(n)\equiv n$ is reached exactly on the family of
permutations constructed by Sturmian words.
\end{abstract}
\end{frontmatter}

\section{Infinite permutations}
Let $S$ be a finite or countable ordered set: we shall consider $S$ equal
either to $\mathbb N$, or to some subset of $\mathbb N$, where
$\mathbb N=\{0,1,2,\cdots\}$.
Let ${\mathcal A}_S$ be the set of all sequences of pairwise distinct reals defined
on $S$. Define an equivalence relation $\sim$ on ${\mathcal A}_S$
as follows: let $a,b \in {\mathcal A}_S$, where $a=\{a_s\}_{s\in S}$ and
$b=\{b_s\}_{s\in S}$; then $a \sim b$ if and only if for all $s,r
\in S$ the inequalities $a_s < a_r$ and $b_s<b_r$ hold or do not
hold simultaneously. An equivalence class from ${\mathcal A}_S / \sim$ is
called an {\it \mbox{($S$-)}permutation}. If an $S$-permutation
$\alpha$ is realized by a sequence of reals $a$, we
denote $\alpha=\overline{a}$. In particular, a
$\{1,\ldots,n\}$-permutation always has a representative with all
values in $\{1,\ldots,n\}$, i.~e., can be identified with a usual
permutation from $S_n$.

In equivalent terms, a permutation can be considered as a linear ordering of $S$
which may differ from the ``natural'' one. That is, for $i,j \in S$, the natural
order between them corresponds to $i<j$ or $i>j$, while the ordering we intend to
define corresponds to $\alpha_i <\alpha_j$ or $\alpha_i >\alpha_j$. We shall also
use the symbols $\gamma_{ij}\in \{<,>\}$ meaning the relations between $\alpha_i$
and $\alpha_j$, so that by definition we have $\alpha_i \gamma_{ij} \alpha_j$ for
all $i \neq j$.

We are interested in properties of infinite permutations analogous to those of infinite words, for example, periodicity and complexity. A permutation $\alpha=\{\alpha_s\}_{s\in S}$ is called {\it
$t$-periodic}
if for all $i,j$ and $n$ such that $i,j,i+nt, j+nt \in S$
we have $\gamma_{ij}=\gamma_{i+nt,j+nt}$. In particular, if $S=\mathbb N$, this definition is equivalent to a more standard one:
a permutation is $t$-periodic if for all $i,j$ we have $\gamma_{ij}=\gamma_{i+t,j+t}$.
A permutation is called {\it ultimately $t$-periodic} if these
equalities hold provided
that $i,j>n_0$ for some $n_0$. This definition is analogous to that for words: an infinite word $w=w_1w_2\cdots$ on an alphabet $\Sigma$ is $t$-periodic if $w_i=w_{i+t}$ for all $i$ and is ultimately $t$-periodic if $w_i=w_{i+t}$ for all $i\geq n_0$ for some $n_0$.

In a previous paper by Fon-Der-Flaass and Frid \cite{ff}, all periodic $\mathbb N$-per\-mu\-ta\-tions have been characterized; in particular, it has been shown that there exists a countable number of distinct $t$-periodic permutations for each $t \geq 2$. For example, for each $n$ the permutation with a representative sequence
\[-1,\; 2n-2,\; 1, \; 2n, \; 3, \; 2n+2, \; 5, \; 2n+4, \; \ldots\]
is 2-periodic, and all such permutations are distinct. So, the situation with periodicity differs from that for words, since the number of distinct $t$-periodic words on a finite alphabet of cardinality $q$ is clearly finite (and is equal to $q^t$).

A set $T=\{0,m_1,\ldots, m_{k-1}\}$ of cardinality $k$, where $0=m_0<m_1<\cdots
<m_{k-1}$, is called a {\em ($k$-)window}. It is natural to define {\em $T$-factors}
of an $S$-permutation $\alpha$ as
restrictions of $\alpha$ to $T+n$, $n \in \mathbb{N}$, considered as permutations
on $T$. Such a projection is well-defined for a given $n$ if and only if $T+n \subseteq S$, and is denoted by
$\alpha_{T+n}=\alpha_n \alpha_{n+m_1}\cdots \alpha_{n+m_{k-1}}$. We call the
number of distinct $T$-factors of $\alpha$ the {\em $T$-complexity} of $\alpha$
and denote it by $p_{\alpha}(T)$.


In particular, if $T=\{0,1,2,\ldots,n-1\}$, then $T$-factors of an $\mathbb N$-permutation $\alpha$ are called just {\em factors} of $\alpha$ and are analogous to factors (or subwords) of infinite words. They are denoted by $\alpha_{[i..i+n)}$ or, equivalently,
$\alpha_{[i..i+n-1]}=\alpha_i \alpha_{i+1}\cdots \alpha_{i+n-1}$, and their number is called the
{\em factor complexity} $f_{\alpha}(n)$ of  $\alpha$. This function is analogous to the subword complexity $f_w(n)$ of infinite words which is equal to the number of different words $w_{[i..i+n)}$ of length $n$ occurring in an infinite word $w$ (see \cite{fer} for a survey). However, not all the properties of these two functions are similar \cite{ff}. Consider in particular the following classical theorem.
\begin{theorem}\label{w}
An infinite word $w$ is ultimately periodic if and only if $f_w(n)=C$ for some constant $C$ and all sufficiently large $n$. If $w$ is not ultimately periodic, then $f_w(n)$ is
increasing and satisfies $f_w(n)\geq n+1$.
\end{theorem}


Only the first statement of Theorem \ref{w} has an analogue for permutations; as for the second one, the situation with permutations is completely different.

\begin{theorem}
\label{per_compl}\cite{ff}
Let $\alpha$ be an $\mathbb N$-permutation; then
$f_{\alpha}(n)\leq C$ if and only if $\alpha$ is
ultimately periodic. At the same time, for each unbounded nondecreasing
 function $g(n)$,
 there exists a $\mathbb N$-permutation $\alpha$ with
$f_{\alpha}(n)\leq g(n)$ for all $n\geq N_0$ which is not ultimately periodic.
\end{theorem}
The supporting example of a permutation with low complexity can be defined by
the inequalities $\alpha_{2n} < \alpha_{2n+2} < \alpha_{2n+1} < \alpha_{2n+3}$ for all $n\geq 0$, and
 $\alpha_{2n_k} < \alpha_{2k+1} < \alpha_{2n_k+2}$
for some sequence $\{n_k\}_{k=0}^{\infty}$ which
grows sufficiently fast.

In this paper we study the properties of another complexity function, namely, {\em maximal pattern complexity}
\[p^*_\alpha(n)=\max_{\#T=n}p_\alpha(T).\]
The analogous function $p^*_w(n)$ for infinite words was defined in 2002 by Kamae and Zamboni \cite{kz} where the following statement was proved:

\begin{theorem}\cite{kz}
An infinite word $w$ is not ultimately periodic if and only if $p^*_w(n)\geq 2n$ for all $n$.
\end{theorem}
Infinite words of maximal pattern complexity $2n$ include rotation words \cite{kz} and also some words built by other techniques \cite{kz2}. The classification of all words of maximal pattern complexity $2n$ is an open problem \cite{k3}.

In this paper, we prove analogous results for infinite permutations and furthermore, prove that in the case of permutations, lowest maximal pattern complexity is achieved only in the precisely described ``Sturmian'' case.

\section{Lowest complexity}

First of all, we prove a lower bound for the maximal pattern complexity of a non-periodic infinite permutation.

\begin{theorem}\label{t_n}
An infinite permutation $\alpha$ is not ultimately periodic if and only if $p^*_{\alpha}(n)\geq n$ for any $n$.
\end{theorem}
{\sc Proof.} Clearly, if a permutation is ultimately periodic, its maximal pattern complexity is ultimately constant, and thus the ``if'' part of the proof is obvious. Now suppose that $p^*_{\alpha}(l)<l$ for some $l$; we shall prove that $\alpha$ is ultimately periodic.

Since $p^*_{\alpha}(1)=1$ (there is exactly one permutation of length one),
and the function $p^*$ is non-decreasing, we see that $p^*_{\alpha}(l)<l$ implies
that $p^*_{\alpha}(n+1)=p^*_{\alpha}(n)$ for some $n \leq l$. Consider an
$n$-window $T=(0,m_1,\ldots,m_{n-1})$ such that $p_{\alpha}(T)=p^*_{\alpha}(n)$;
the equality $p^*_{\alpha}(n+1)=p_{\alpha}(T)$ implies that for each
$T'=(0,m_1,\ldots,m_{n-1}, m_n)$ with $m_n >m_{n-1}$ we have  $p_{\alpha}(T)=p_{\alpha}(T')$,
that is, each $T$-permutation can be extended to a $T'$-permutation in a unique way.
Clearly, there exist two equal factors of length $2m_{n-1}$ in $\alpha$: say,
$$\alpha_{[k..k+2m_{n-1})}=\alpha_{[k+t..k+t+2m_{n-1})}$$
for some positive $t$ and non-negative $k$.
 We shall prove that $\alpha$ is ultimately $t$-periodic, namely, that $\gamma_{ij}=\gamma_{i+t,j+t}$ for all $i,j$ with $k\leq i< j$. The proof will use the induction on the pair $i,j$ starting by the pairs $i,j$ with $k\leq i<j< k+2m_{n-1}$, for which our statement holds since $\alpha_{[k..k+2m_{n-1})}=\alpha_{[k+t..k+t+2m_{n-1})}$.

Now for the induction step: for some $M\ge 2m_{n-1}$, suppose that $\gamma_{ij}=\gamma_{i+t,j+t}$ for all $k\leq i<j<k+M$, that is, $\alpha_{[k..k+M)}=\alpha_{[k+t..k+t+M)}$. We are going to prove that $\gamma_{i,k+M}=\gamma_{i+t, k+t+M}$ for all $i\in\{k,\ldots, k+M-1\}$, and thus  $\alpha_{[k..k+M+1)}=\alpha_{[k+t..k+t+M+1)}$.

Indeed, condider the case $i\in \{k,\ldots,k+M-m_{n-1}-1\}$ first. Then $\alpha_{T+i}$ is a $T$-factor of $\alpha_{[k..k+M)}$ and $\alpha_{T+i+t}$ is a $T$-factor of $\alpha_{[k+t..k+t+M)}$ standing at the same position. So, these $T$-factors of $\alpha$ are equal, and due to the choice of $T$, so are their extensions $\alpha_{T'+i}$ and $\alpha_{T'+i+t}$, where $T'=(0,m_1,\ldots, m_{n-1},M-i)$. In particular, the first and last elements of $\alpha_{T'+i}$ and $\alpha_{T'+i+t}$ are in the same relationship: $\gamma_{i,k+M}=\gamma_{i+t,k+t+M}$, which is what we needed.

Now if $i\in \{k+M-m_{n-1},\ldots,k+M-1\}$, we consider $\alpha_{T+i-m_{n-1}}$ which is a $T$-factor of $\alpha_{[k..k+M)}$ with the last element $\alpha_i$, and  $\alpha_{T+i+t-m_{n-1}}$ which is a $T$-factor of $\alpha_{[k+t..k+t+M)}$ with the last element $\alpha_{i+t}$. They are equal, and so are their extensions $\alpha_{T'+i-m_{n-1}}$ and $\alpha_{T'+i+t-m_{n-1}}$, where $T'=(0,m_1,\ldots, m_{n-1},M-i+m_{n-1})$. In particular, the next to last and the last elements of these $T$-permutations are in the same  relationship:  $\gamma_{i,k+M}=\gamma_{i+t,k+t+M}$.

So, $\gamma_{i,k+M}=\gamma_{i+t,k+t+M}$ for all $i \in \{k,\ldots, k+M-1\}$; together with the induction hypothesis it means that $\alpha_{[k..k+M+1)}=\alpha_{[k+t..k+t+M+1)}$. Repeating the induction step we get that $\gamma_{ij}=\gamma_{i+t,j+t}$ for all $k\leq i<j$, that is, the permutation $\alpha$ is ultimately $t$-periodic.  \hfill $\Box$

\section{Sturmian permutations}\label{sec_def}
A one-sided infinite word $w=w_0w_1w_2\cdots$ on the alphabet $\{0,1\}$ is called
{\em Sturmian} if its subword complexity $f_w(n)$ is equal to $n+1$ for all $n$.
Sturmian words have a number of equivalent definitions \cite{sturm}; we
shall need two more of them. First, Sturmian words are exactly aperiodic
{\em balanced} words which means that for each length $n$, the number of $1$'s in
factors of $w$ of length $n$ takes only two successive values. Second, Sturmian
words are exactly irrational {\em mechanical} words which means that there exists
some irrational $\sigma \in (0,1)$ and some $\rho \in [0,1)$ such that for all $i$ we have
\begin{align} \label{st}
w_i=\lfloor \sigma (i+1) + \rho \rfloor - \lfloor \sigma i + \rho \rfloor & \mbox{~or} \\
w_i=\lceil \sigma (i+1) + \rho \rceil - \lceil \sigma i + \rho \rceil &.
\end{align}
These definitions coincide if $\sigma i + \rho$ is never integer;
if it is integer for some (unique) $i$, the sequences built by these two formulas differ
in at most two successive positions. So, we distinguish {\em lower} and
{\em upper} Sturmian words according to the choice of $\lfloor \cdot \rfloor$ or
$\lceil \cdot \rceil$ in the definition. A word on any other binary alphabet is
called Sturmian if it is obtained from a Sturmian word on $\{0,1\}$ by renaming
symbols. Here $\sigma$ is called the {\em slope} of the word $w$.

Now let us define a {\em Sturmian permutation} $\alpha(w,x,y)=\alpha=\overline{a}$
associated with a Sturmian word $w$ and positive numbers $x$ and $y$ by its
representative sequence $a$, where $a_0$ is a real number and for all
$i\geq 0$ we have
\[a_{i+1}=\begin{cases} a_i+x, \mbox{~if~} w_i=0,\\
a_i-y, \mbox{~if~} w_i=1. \end{cases}\]
Clearly, such a permutation is well-defined if and only if we never have
$kx\neq ly$ if $k$ is the number of $0$'s and $l$ is the number of $1$'s in some factor of $w$; and in particular if $x$ and $y$ are rationally independent.

Note that a factor of $w$ of length $n$ corresponds to a factor of $\alpha$ of length $n+1$, and the correspondence is one-to-one. So, we have $f_{\alpha}(n)=n$ for all $n$. In fact, we are going to prove that the maximal pattern complexity of $\alpha$ is also equal to $n$, and thus the lower bound in Theorem \ref{t_n} is precise.

\begin{theorem}
For each Sturmian permutation $\alpha$ we have $p^*_{\alpha}(n)\equiv n$.
\end{theorem}
{\sc Proof.} Let us start with the situation when $x=\sigma$ and $y=1-\sigma$. This case has been proved by M. Makarov in \cite{mak2}, but we give a proof here for the sake of completeness.

If we take $a_0=\rho$, then by the definition of the Sturmian word,
$a_i=\{\sigma i +\rho\}$ holds in the case that $w$ is a lower Sturmian word,
and $a_i=1-\{1-\sigma i -\rho\}$ holds in the case that $w$ is an upper Sturmian
word. Here $\{x\}$ stands for the fractional part of $x$.
In what follows, we consider lower Sturmian words without loss of generality.

Consider a $k$-window $T=\{0,m_1,\ldots,m_{k-1}\}$ and the set of $T$-factors
$\alpha_{T+n}=\overline{\{\sigma n +\rho\},\{\sigma (n+m_1) +\rho\},\ldots,\{\sigma (n+m_{k-1}) +\rho\}}$ for all $n$. Since the set of $\{\sigma n +\rho\}$ is dense in $[0,1]$, the set of $T$-factors is equal to the set of all permutations $\overline{t,\{t+\sigma m_1\},\ldots,\{t+\sigma m_{k-1}\}}$ with $t \in [0,1]$

Let us arrange the points $\{t+\sigma m_i\}$ ($i=0,\ldots,k-1$) on the unit circle, that is the interval $[0,1]$ with the points 0 and 1 identified (recall that $m_0=0$ by definition). Then, the arrangement partitions the unit circle into $k$ arcs. Since the arrangements for different $t$'s are different only by rotations, the permutation defined by the points is determined by indicating the arc which contains $0=1$. Since there exist $k$ arcs, there are exactly $k$  different permutations defined by the points $\{t+\sigma m_i\}$ ($i=0,\ldots,k-1$) with different $t$'s. Thus, $p_{\alpha}(T)=k$. Since the window $T$ was chosen arbitrarily, we have $p^*_{\alpha}(k)=k$.

Now consider the general case of arbitrary $x$ and $y$.
Let us keep the notation $\gamma_{ij}$ for the relation between $\alpha(w,\sigma,1-\sigma)_i$ and $\alpha(w,\sigma,1-\sigma)_j$, and denote the relation between $\alpha(w,x,y)_i$ and $\alpha(w,x,y)_j$ by $\delta_{ij}$.

Recall that the {\em weight} of a binary word $u$ is the number of 1's in it, denoted by $|u|_1$.
 By the definition of $\alpha$, we have $\delta_{i,i+n}=\delta_{j,j+n}$ if $w_{[i..i+n)}$ and $w_{[j..j+n)}$ have the same weight. Note also that the weight of a factor of $w$ of length $n$ is either equal
 to $\lfloor n \sigma \rfloor$ or to $\lceil n \sigma \rceil$. In $\alpha(w,\sigma,1-\sigma)$, 
the converse also holds: words $w_{[i..i+n)}$ and $w_{[j..j+n)}$  of the same length $n$ but of different weight always correspond to $\gamma_{i,i+n}\neq \gamma_{j,j+n}$, since $(n-\lfloor n \sigma \rfloor)\sigma - \lfloor n \sigma \rfloor (1-\sigma)= n\sigma - \lfloor n \sigma \rfloor>0$ and $(n-\lceil n \sigma \rceil)\sigma - \lceil n \sigma \rceil (1-\sigma)= n\sigma - \lceil n \sigma \rceil < 0$. In the general case, words of different weights may correspond to the same relation. But anyway for all $i$, $j$, and $n$ the equality $\gamma_{i,i+n}=\gamma_{j,j+n}$ implies that $\delta_{i,i+n}=\delta_{j,j+n}$. Thus, for any $k$-window $T$ we see that 
$\alpha(w,\sigma,1-\sigma)_{T+i}=\alpha(w,\sigma,1-\sigma)_{T+j}$ implies $\alpha(w,x,y)_{T+i}=\alpha(w,x,y)_{T+j}$. So, we have
$p_{\alpha(w,x,y)}(T)\leq p_{\alpha(w,\sigma,1-\sigma)}(T)$ and thus
$p^*_{\alpha(w,x,y)}(k)\leq p^*_{\alpha(w,\sigma,1-\sigma)}(k)=k$; at the same
time, $p^*_{\alpha(w,x,y)}(k)\geq k$ since this permutation is not ultimately
periodic. So, $p^*_{\alpha(w,x,y)}(k)= k$, and the theorem is proved.
\hfill $\Box$.

In fact, Sturmian permutations are the only $\mathbb N$-permutations of maximal pattern complexity $n$. In the remaining part of the paper, we are going to prove it. 

\section{Rotation words}
In what follows, we several times use the fact that Sturmian words
form a particular case of so-called {\em rotation words}. Let us
describe them.

Consider the interval $\mathcal{C}=[0,1)$ as a unit circle, which means that we identify
its ends and consider it as the quotient group $\mathbb R / \mathbb
Z$. When working with this group, we consider real numbers modulo
one and write $x \pmod 1$ or just $x$ as well as the fractional part
$\{x\}$.

An interval $I=[x,y)$ on
$\mathcal{C}$ is defined as usual if $0\leq x< y <1$ and as
$\mathcal{C}\backslash [y,x)$ if 
$0\leq y <x <1$. Intervals with other
combinations of parentheses are defined analogously.

Now consider a partition of $\mathcal{C}$ into a finite number of
disjoint intervals $J_0, J_1,\ldots, J_k$, $\cup_{j=0}^{k} J_j =
\mathcal{C}$. Associate with each interval $J_j$ a symbol $a_j$ from a finite
alphabet $A$ (symbols for different intervals may coincide). Let
$I_{a}$ denote the union of intervals corresponding to the symbol $a$. 

Consider a sequence $(x_i)_{i=0}^{\infty}$, $x_i \in \mathcal C$, given by
$x_{i+1}=x_i + \xi \pmod 1$ for some fixed $\xi$, and define an infinite word
$v=v_0\cdots v_n\cdots$ on the alphabet $A$ by $v_i=a \Longleftrightarrow x_i \in I_a$. This word is called a {\em rotation word} on $A$ with the {\em slope} $\xi$ and the initial point $x_0$ induced by the given partition of $\mathcal C$.

Thus, a Sturmian word defined by \eqref{st} is a rotation word induced by a partition of
$\mathcal C$ into the intervals $[0,\sigma)$ and $[\sigma,0)$ (for a
lower Sturmian word; for the upper Sturmian word, the parentheses
are $(\cdot,\cdot]$); with the initial point $x_0=\sigma+\rho$.
Equivalently, we can define it by the partition into the intervals
$[-\sigma-\rho,-\rho)$ and $[-\rho,-\sigma-\rho)$ with the initial
point $0$.

\section{Proof of uniqueness: first step}
Now we shall prove that the described Sturmian permutations are the only
permutations of maximal pattern complexity $p^*_{\alpha}(n)= n$. In the proof,
we shall widely use the table of values $\gamma_{ij}\in\{<,>\}$ of a candidate permutation;
for the sake of convenience, we denote the strings of that table by
$\gamma_i=\gamma_{0,i}\gamma_{1,i+1}\cdots\gamma_{n,i+n}\cdots$ and the
arithmetical subsequences of those strings by
$$\gamma_i^j=\gamma_{j,i+j}\gamma_{i+j,2i+j}\cdots\gamma_{ni+j,(n+1)i+j}\cdots$$
for all $i \in {\mathbb N}$ and $j \in \{0,\ldots,i-1\}$.
Thus, a string $\gamma_i$ consists of elements of $i$ disjoint sequences
$\gamma_i^j$, 
each of them representing the relations between successive 
elements of the permutation $\{\alpha_{ni+j}\}_{n=0}^\infty$.

So, each $\gamma_i^j$ is an infinite word on the alphabet $\{<,>\}$.
We also denote the subword $\gamma_{n,n+1}\gamma_{n+1,n+2} \cdots
\gamma_{n+i-1,n+i}$ by $\gamma_{[n..n+i]}$.

\begin{lemma}\label{l1}
If  $\alpha$ is an infinite permutation with $p^*_{\alpha}(n)\equiv n$, then for all $i>0$ and $j\in\{0,\ldots, i-1\}$ the sequence $\gamma_i^j$ is either 
ultimately periodic or Sturmian.
\end{lemma}
{\sc Proof.} Let us fix some $i$. If $p^*_{\alpha}(n)\equiv n$, then in particular $p_\alpha(T_n)\leq n$, where $T_n=(0,i,2i,\ldots,(n-1)i)$. Thus, the number of different values $\alpha_{j+ik+T_n}$ for different $k$'s are at most $n$, and since the factor $$\gamma_{j+ki,j+(k+1)i} \gamma_{j+(k+1)i,j+(k+2)i}\cdots \gamma_{j+(k+n-2)i,j+(k+n-1)i}$$ of $\gamma_i^j$ contains just a part of information contained in $\alpha_{j+ik+T_n}$, the number of such factors of length $n-1$ is at most $n$ for all $n$. Since the only non-periodic words satisfying this are Sturmian words, the lemma is proved. \hfill $\Box$

In particular, this lemma is valid for $\gamma_1=\gamma_1^0$. In what follows we consider the cases when  $\gamma_1$ is periodic and when it is Sturmian separately.

\section{Proof of uniqueness: Sturmian case}\label{sss}
In this section we assume that the first string $\gamma_1$ of the array $\{\gamma_i^j\}$, describing the relations between successive elements of a permutation $\alpha$ with $p^*_{\alpha}(n)=n$, is a Sturmian word on the alphabet $\{<,>\}$. Let us see what all the other substrings $\gamma_i^j$ are.

We say that an infinite word on $\{<,>\}$ is {\em increasing} (or {\em decreasing}) if it is equal to $<^{\omega}$ (or $>^{\omega}$, respectively). It is called {\em monotonic} if it is either increasing or decreasing. We put ``{\em ultimately}" if it holds after some point. 

\begin{claim}\label{cll}
For each $i>0$ and $j \in \{0,\ldots,i-1\}$ the sequence $\gamma_i^j$ is either
Sturmian or ultimately monotonic.
\end{claim}
{\sc Proof.} Due to Lemma \ref{l1}, it is sufficient to prove that $\gamma_i^j$ cannot be ultimately periodic with the minimal period $t$
 greater than one, that is, we cannot have for any $t$, $m_1$, and $m_2$ that
 $(\gamma_i^j)_{m_1+nt}=\gamma_{j+(m_1+nt)i,j+(m_1+nt+1)i}=<$ for all sufficiently large $n$ and $(\gamma_i^j)_{m_2+nt}=\gamma_{j+(m_2+nt)i,j+(m_2+nt+1)i}=>$ for all sufficiently large $n$. To the contrary, let us suppose this and consider the pattern $T=(0,i,i+1)$. Consider the $T$-permutations $\alpha_{k+T}$ for all $k$. Each of them is determined by the three values:    $\gamma_{k,k+i}$, $\gamma_{k,k+i+1}$, and $\gamma_{k+i,k+i+1}$. Consider first $k=j+(m_1+nt)i$ for all sufficiently large $n$. We see that $\gamma_{k,k+i}$ in this case is equal to $<$, but $\gamma_{k+i,k+i+1}$ takes both values for different $n$'s since the sequence $\gamma_1$ is Sturmian and thus any infinite arithmetic progression in it contains both symbols. Analogously, if $k=j+(m_2+nt)i$, then $\gamma_{k,k+i}$ is ultimately equal to $>$ and $\gamma_{k+i,k+i+1}$ takes both values. So, $T$-permutations $\alpha_{k+T}$ take at least four values, which means that $p^*_{\alpha}(3)\geq 4$, contradicting to the assumption that $p^*_{\alpha}(n)=n$. \hfill $\Box$

\begin{claim}
Given $i$, if $\gamma_i^{j_1}$ is Sturmian for some $j_1$, then $\gamma_i^j$ is Sturmian for any $j=0,\ldots,i-1$.
\end{claim}
{\sc Proof.} Due to the previous claim, the opposite would mean that some of $\gamma_i^j$ were ultimately monotonic. Suppose without loss of generality that 
$\gamma_i^j$ is ultimately increasing, and let $n$ be the greatest number of successive symbols $<$ in $\gamma_i^{j_1}$ (clearly it is finite). Consider the pattern $T_{n+2}=(0,i,\ldots,ni,(n+1)i)$ of length $n+2$. For different $k$ equal to $j_1$ modulo $i$, the number of different $\alpha_{k+T_{n+2}}$'s is at least $n+2$ since $\gamma_i^{j_1}$, the sequence describing the relations between the successive elements of $T_{n+2}$, is Sturmian. Moreover, since $<^{n+1}$ is not contained in $\gamma_i^{j_1}$, while it is contained in $\gamma_i^{j}$, $\alpha_{k+T_{n+2}}$ can take at least $n+3$ different values, contradicting to the assumption that $p^*_{\alpha}(k)=k$. \hfill $\Box$

\begin{claim}
Suppose that $\gamma_i^{j_1}$ is ultimately increasing (ultimately decreasing) for some $j_1$. Then, $\gamma_i^j$ is ultimately increasing (ultimately decreasing, respectively) for any $j=0,1,\ldots,i-1$.
\end{claim}
{\sc Proof.} Due to the previous claims, the opposite would mean exactly that 
$\gamma_i^{j_2}$ is ultimately decreasing for some $j_2 \in\{0,\ldots,i-1\}$, while $\gamma_i^{j_1}$ is ultimately increasing. Now consider once again the pattern $T=(0,i,i+1)$ and like in Claim \ref{cll} observe that the pair  $(\gamma_{n,n+i},\gamma_{n+i,n+i+1})$ which contains a part of information of $\alpha_{n+T}$, takes at least two different values $(<,<)$ and $(<,>)$  when $n=j_1~(mod~i)$. Also, it takes two values $(>,<)$ and $(>,>)$ when $n=j_2~(mod~i)$. So, $p^*_{\alpha}(3)\geq p_{\alpha}(T)\geq 4$, a contradiction. \hfill $\Box$

\begin{claim}
If $\gamma_i$ and $\gamma_j$ are ultimately monotonic, then they are ultimately increasing or ultimately decreasing, simultaneously.
\end{claim}
{\sc Proof.} Suppose the opposite: say, 
$\gamma_i$ is ultimately increasing and $\gamma_j$ is ultimately decreasing. It means that for a sufficiently large $k$ we have 
$\alpha_k<\alpha_{k+i}<\alpha_{k+2i}<\ldots <\alpha_{k+ji}$, and at the same time, $\alpha_k>\alpha_{k+j}>\alpha_{k+2j}>\ldots >\alpha_{k+ij}$, a contradiction. \hfill $\Box$

Therefore, the set of positive integers is divided into two classes $S$ and $M$:
a number $i$ belongs to $S$ if all $\gamma_i^j$ are Sturmian, and to $M$ if
$\gamma_i$ is ultimately monotonic. Due to the previous claim, all $\gamma_i$ with $i \in M$ are ultimately increasing or ultimately decreasng, simultaneously, and without loss of generality we may assume that they are ultimately decreasing. Now let us specify what kind of Sturmian words $\gamma_i^j$\,'s are. 

Let the slope of the Sturmian word $\gamma_1$ be equal to $\sigma$ and the
initial point be $\rho$.
Without loss of generality we assume that the word is
lower Sturmian: this means precisely that
\[\gamma_{n,n+1}=\begin{cases} < , \mbox{~if~} \{\sigma (n+1)
+\rho\}<\sigma,
\\ > , \mbox{~otherwise~}. \end{cases}\]
In other words, $\gamma_{n,n+1}=<$ if and only if $\sigma n \in
[-\sigma-\rho, -\rho) \bmod 1$. Moreover, $\gamma_{n+1,n+2}=<$ if
and only if $\sigma n \in
[-2\sigma-\rho, -\sigma-\rho) \bmod 1$, etc.: we see that the word
$\gamma_{[n..n+i]}$ is determined by the position of the point
$\sigma n \in {\mathcal C}$ with respect to the points $-\rho,
-\sigma-\rho,\ldots,-i\sigma-\rho$.

Let us fix some $i$. We know that
\[\# \{\alpha_{[n..n+i]}|n \in {\mathbb N}\} = f_{\alpha}(i+1) \leq
p^*_\alpha(i+1)=i+1.\]
On the other hand, we have 
$$
\# \{\alpha_{[n..n+i]}|n \in {\mathbb N}\}\geq 
\# \{\gamma_{[n..n+i]}|n \in {\mathbb N}\}=i+1. 
$$
Hence, we have $\# \{\alpha_{[n..n+i]}|n \in {\mathbb N}\}=
\# \{\gamma_{[n..n+i]}|n \in {\mathbb N}\}$. 
It follows that the whole 
permutation $\alpha_{[n..n+i]}$, and in particular the relation
$\gamma_{n,n+i}$, is uniquely determined by $\gamma_{[n..n+i]}$
and thus by the position of the point 
$\sigma n \in {\mathcal C}$ with respect to the points $-\rho,
-\sigma-\rho,\ldots,-i\sigma-\rho$ modulo 1.

For $i \in M$ this implies that the sequence $\gamma_i$ is monotinic, not only ultimately monotonic.

For $i \in S$ this means that $\gamma_i$ is a rotation word on $\{<,>\}$ with the slope
$\sigma$ starting at 0, and the partition of $\mathcal{C}$ by 
the set of intervals of type $[~,~)$ bounded by the points $-\rho,
-\sigma-\rho,\ldots,-i\sigma-\rho$ modulo 1. And for each $j=0,\ldots,i-1$, the
word $\gamma_i^j$ is a rotation word on $\{<,>\}$ corresponding to the same partition of $\mathbb{C}$ by the intervals, with the slope $i\sigma$ starting at $j\sigma$. Here, we have not yet specified the correspondance between intervals and $\{<,>\}$.

\begin{claim}\label{c5}
Assume that $i \in S$, which means that $\gamma_i^j$ are Sturmian words for all $j$. Let $I_<$ be the union of the above intervals corresponding to $<$. Then, $I_<$ is an interval in $\mathcal{C}$ of length $\{i \sigma\}$ or $1-\{i \sigma\}$.
\end{claim}
{\sc Proof.} Note that $\gamma_i^j$ is a rotation word on $\{<,>\}$ with the slope $i\sigma$ starting at $j\sigma$ corresponding to the partition defined by the set of points $S=\{\{-\rho\},\{-\sigma-\rho\},\ldots,\{-i\sigma-\rho\}\}$. Note that $S\cap(S+i\sigma)=\{\{-\rho\}\}$ and $S\cap(S+ki\sigma)=\emptyset$ 
for any $k=2,3,\ldots$.  

Suppose that the conclusion in the Claim does not hold. Then, there exist $u,v\in S$ with $u<v$ in the boundary of $I_<$ such that $v-u\ne \{i\sigma\}$ and  $v-u\ne 1-\{i\sigma\}$. Let $\mathbb{P}=\{I_<,I_>\}$ be the partition of $\mathcal{C}$ and  
$$\mathbb{P}_{k+1}=\mathbb{P}\vee(\mathbb{P}-i\theta)\vee\ldots\vee(\mathbb{P}-ki\theta)$$
be the refinement of partitions. There exists $k>0$ such that $u$ and $v$ are in the interiors of distinct elements of the partition 
$$\mathbb{P}'=(\mathbb{P}-i\theta)\vee\ldots\vee(\mathbb{P}-ki\theta).$$
Then, we have $\#\mathbb{P}_{k+1}\ge\#\mathbb{P}'+2$ since each of the points $u$ and $v$ increases $\#\mathbb{P}_{k+1}$ from $\#\mathbb{P}'$ by 1. On the other hand, $\#\mathbb{P}'=\#\mathbb{P}_k=k+1$ and $\#\mathbb{P}_{k+1}=k+2$ hold since $\gamma_i^j$ is a Sturmian word. Thus, we have a contradiction. 
\hfill $\Box$

By Claim \ref{c5}, only 2 cases are possible. That is, either 
$I_<=[\{-\rho\},\{-i\sigma-\rho\})$ and $I_>=[\{-i\sigma-\rho\},\{-\rho\})$ 
or $I_>=[\{-\rho\},\{-i\sigma-\rho\})$ and $I_<=[\{-i\sigma-\rho\},\{-\rho\})$. Hence, there are only two Sturmian
words on $\{<,>\}$, satisfying our properties, and they are 
obtained from the other by exchanging the symbols.

In fact, we can describe the Sturmian words obtained in the above 
in a more direct way. To do it, for each $j=0,\ldots,i-1$ consider 
the word $v^j=v_0\cdots v_n \cdots$ defined by
$v_n=|\gamma_{[j+ni..j+(n+1)i]}|_<$, where $|w|_a$ denotes the number
of occurrences of a symbol $a$ in the word $w$. As it follows from the definition of $\gamma_1$, the word
$v$ is binary on the alphabet $\{q_i,q_i+1\}$, where $q_i=\lfloor
\sigma i \rfloor$. It is not periodic since $\gamma_1$ is Sturmian.
Moreover, the word $v$ is balanced since
$|\gamma_{[j+ni..j+k(n+1)i]}|_<=q_ik+|v_{[n..n+k]}|_{q_i+1}$ also takes
only two values for a fixed $k$, and so does $|v_{[n..n+k]}|_{q_i+1}$.
But non-periodic balanced words are exactly Sturmian words. So, $v$
is Sturmian, and its symbol $v_n$ is determined by
$\gamma_{[j+ni..j+(n+1)i]}$. As we have shown in the previous
paragraph, it means that $\gamma_i^j$ is obtained from $v$ by
renaming symbols, that is, each its symbol $\gamma_{j+ni,j+(n+1)i}$
is determined by the number of symbols $<$ in $\gamma_{[j+ni..j+(n+1)i]}$ 
independently of $j$, which is either $q_i$ or $q_i+1$. 

Thus, there is a mapping, say $\rho_i$, from $\{q_i,q_i+1\}$ to $\{<,>\}$ 
such that $\rho_i(|\gamma_{[m..m+i]}|_<)=\gamma_{m,m+i}$ 
for any $i,m$. 

Consider the case $i\in M$. Since $\gamma_{m,m+i}$ is independent of $m$ 
for any large $m$, we have $\rho_i(q_i)=\rho_i(q_i+1)$. Thus, $\rho_i$ takes only one value and it holds that $\gamma_i$ is not only ultimately monotonic, but also monotonic for any $i\in M$. 

Consider the case $i\in S$. 
Since $i\in S$, and thus $\gamma_{m,m+i}$ can take both values $<$ and $>$, $\rho_i$ is a bijection. 
Suppose first that $\rho_i(q_i)=<$. Consider a factor
$\gamma_{[m..m+i+1]}$ of $\gamma_1$ starting with $>$ and
ending with $<$, so that $|\gamma_{[m..m+i]}|_<=q_i$ and
$|\gamma_{[m+1..m+i+1]}|_<=q_i+1$. We have $\gamma_{m,m+i}=<$, and
thus $\alpha_{m+1}<\alpha_m<\alpha_{m+i}<\alpha_{m+i+1}$. At the
same time, $\gamma_{m+1,m+i+1}=>$, that is,
$\alpha_{m+1}>\alpha_{m+i+1}$. A contradiction to our assumption. 
Hence, we have $\rho_i(q_i)=>$ and $\rho_i(q_i+1)=<$

At last, note that $|\gamma_{[m..m+i]}|_<=q_i$ if and only if $\lfloor
\sigma(m+i)+\rho \rfloor-\lfloor \sigma m + \rho\rfloor = \lfloor
\sigma i\rfloor$, which is equivalent to the inequality $\{\sigma m +
\rho\}< \{\sigma(m+i)+\rho\}$.

We have proved

\begin{claim}
For $i\in M$, $\gamma_i$ is monotonic. 
For $i\in S$, we have $\gamma_{m,m+i}=<$ if and only if
$|\gamma_{[m..m+i]}|_<=q_i+1$, that is, if and only if
$\{\sigma(m+i)+\rho\}<\{\sigma m + \rho\}$. 
\end{claim}

Taken together, the claims above mean that a permutation $\alpha$ of maximal
pattern complexity $p^*_{\alpha}(n)=n$, such that  the upper raw
$\gamma_1$ is a Sturmian word, is uniquely determined by
\begin{itemize}
\item
the Sturmian word $\gamma_1$, and in particular its parameters $\sigma$ and
$\rho$;
\item
the partition of $\mathbb N$ into $S$ and $M$;
\item
the type of (all the words) $\gamma_i$ with $i\in M$: in what
follows we assume without loss of generality that they all are decreasing.
\end{itemize}

However, it is not difficult to see that given a word $\gamma_1$, we
cannot choose the partition $\mathbb N = S \cup M$ arbitrarily. Let
us consider restrictions which we must put on it.

Suppose first that $i,j \in M$. It means that for all large $k$ we
have $\alpha_k>\alpha_{k+i}>\alpha_{k+i+j}$. Since a linear
order is always transitive, this means that
$\alpha_k>\alpha_{k+i+j}$ and thus $i+j \in M$, giving us the
following condition:

\begin{equation}\label{c1}
i, j \in M \Longrightarrow i+j \in M.
\end{equation}

To state other conditions, let us return to the number $q_i=\lfloor i \sigma \rfloor$. Recall that the number of symbols $<$ in the factors of $\gamma_1$ of length $i$ is either $q_i$ or $q_i+1$. 
Since 
$$q_{i+j}+\{(i+j)\sigma\}=(i+j)\sigma=i\sigma+j\sigma=q_i+\{i\sigma\}+q_j+\{j\sigma\},$$
we have $q_{i+j}-q_i-q_j=\{i\sigma\}+\{j\sigma \}-\{(i+j)\sigma\}$. 
Hence, $q_{i+j}-q_i-q_j>0$ if and only if 
$\{i\sigma\}+\{j\sigma \}-\{(i+j)\sigma\}>0$. 
The former is equivalent to $q_{i+j}=q_i+q_j+1$ and the latter is equivalent to $\{i\sigma\}+\{j\sigma \}>1$. Thus, we have 
$$q_{i+j}=q_i+q_j+1\mbox{ if and only if }\{i\sigma\}+\{j\sigma \}>1.$$

Assume that $i,j\in S$ and $\{i\sigma\}+\{j\sigma \}>1$. There exist infinitely many $k$'s such that $|\gamma_{[k..k+i+j]}|_<=q_{i+j}+1$ since $\gamma_1$ is a Sturmian word. On the other hand, we have $q_{i+j}=q_i+q_j+1$ since $\{i\sigma\}+\{j\sigma \}>1$. It follows that $|\gamma_{[k..k+i]}|_<=q_i+1$ and $|\gamma_{[k+i..k+i+j]}|_<=q_j+1$ since $|\gamma_{[k..k+i]}|_<\leq q_i+1$, $|\gamma_{[k+i..k+i+j]}|_<\leq q_j+1$ and $|\gamma_{[k..k+i+j]}|_<=q_{i+j}+1=q_i+1+q_j+1$. Since $i,j\in S$, this implies that $\alpha_k<\alpha_{k+i}<\alpha_{k+i+j}$ and $i+j$ cannot be in $M$. Hence, $i+j\in S$. 
\begin{equation}\label{c2}
i, j \in S \mbox{~and~} \{i \sigma\} + \{j \sigma \} > 1 \Longrightarrow i+j \in S.
\end{equation}

Now consider the situation when $i+j \in S$ and a word of length $i+j$ in $\gamma_1$ with $q_{i+j}+1$ occurrences of $<$ ends by a suffix of length $j$ with only $q_j$ occurrences of $<$. This is possible if and only if $q_{i+j}=q_{i}+q_j$, that is, $\{i\sigma\}+\{j\sigma\}=\{(i+j)\sigma\}<1$. There exists $k$ such that $|\gamma_{[k..k+i+j]}|_<=q_{i+j}+1$. Then we have $\alpha_k<\alpha_{k+i+j}$ since $i+j \in S$ and $\alpha_{k+i}>\alpha_{k+i+j}$ since $|\gamma_{[k+i..k+i+j]}|_<=q_j$: here it does not matter if $j\in S$ or $j \in M$. Thus, by transitivity $\alpha_k<\alpha_{k+i}$ holds, which means in particular that $i \in S$. We have proved that
\begin{equation}\label{c3}
i+ j \in S \mbox{~and~} \{i \sigma\} + \{j \sigma \} < 1 \Longrightarrow i \in S.
\end{equation}

Note that $i$ and $j$ in this condition are treated symmetrically, so in fact, $j$ also belongs to $S$.

Now using the conditions \ref{c1}--\ref{c3} we can prove 

\begin{claim}\label{dd}
For each $s \in S$ and $m \in M$ we have 
\[\frac{1-\{ m\sigma \}}{m}< \frac{1-\{ s\sigma \}}{s}.\]
\end{claim}
{\sc Proof.} Suppose to the contrary that 
\begin{equation}\label{contr}
\frac{1-\{ m\sigma \}}{m}\geq \frac{1-\{ s\sigma \}}{s}
\end{equation}
 for some $s \in S$ and $m \in M$, and choose a minimal counter-example, so that the sum of $s$ and $m$ is the least possible.

Suppose first that $s>m$. Then $s-m \in S$ due to \eqref{c1}. Moreover, since $m \in M$, we do not get into Condition \eqref{c3}, and thus $\{m \sigma\} + \{(s-m) \sigma \} > 1$, that is, $\{m \sigma\} + \{(s-m) \sigma \}=\{s \sigma\} +1$. It can be checked directly using \eqref{contr} that $\displaystyle \frac{1-\{(s-m)\sigma\}}{s-m}=\frac{\{m \sigma\}-\{s\sigma\}}{s-m}\leq \frac{1-\{m \sigma\}}{m}$, so that $s-m \in S$ and $m \in M$ form a counter-example less than the initial one, contradicting to its minimality.

Now suppose that $m>s$. Then \eqref{contr} immediately implies that $\{m \sigma\}< \{s \sigma\}$ (equality being impossible since $\sigma$ is irrational), and thus $\{s \sigma\} + \{(m-s) \sigma \}=\{m \sigma\}+1$ (not $\{m \sigma\}$). Due to \eqref{c2}, we have $m-s \in M$ since otherwise we would have $m \in S$. Now we again can see that $s$ and $m-s$ give a counter-example less than the initial one since $\displaystyle \frac{1-\{(m-s)\sigma\}}{m-s}=\frac{\{s \sigma\}-\{m \sigma\}}{m-s}\geq\frac{1-\{s \sigma\}}{s}$ due to \eqref{contr}. \hfill $\Box$.

Now note that $\displaystyle \frac{1-\{i \sigma\}}{i} \to 0 $ with $i \to \infty$. Note also that the set $S$ is not empty since $1 \in S$. So, Claim \ref{dd} means that either $S=\mathbb N$, or there exists some $d\in (0,1)$ such that $i \in S$ if and only if $\frac{1-\{ i\sigma \}}{i}>d$, and $i\in M$ if and only if $\frac{1-\{ i\sigma \}}{i}<d$. This parameter $d$ together with the word $\gamma_1$ and the fact that the monotonic strings of the table $\gamma$ are decreasing, completely defines the permutation $\alpha$. Note that the situation when  $S=\mathbb N$ just corresponds to $d=0$.

It remains to check that $\alpha=\alpha(\gamma_1,1-\sigma-d,\sigma+d)$. Here, we just treat each symbol $<$ in $\gamma_1$ as $0$ and $>$ as $1$ to use the definition of a Sturmian permutation from Section \ref{sec_def}. Indeed, $\gamma_{k,k+i}=>$ if and only if $|\gamma_{[k..k+i]}|_<=q_i$ which implies $\alpha(\gamma_1,1-\sigma-d,\sigma+d)_k>\alpha(\gamma_1,1-\sigma-d,\sigma+d)_{k+i}$ since $q_i(1-\sigma-d)-(i-q_i)(\sigma+d)<0$. On the other hand, $\gamma_{k,k+i}=<$ if and only if $|\gamma_{[k..k+i]}|_<=q_i+1$ and $i\in S$ which implies $\alpha(\gamma_1,1-\sigma-d,\sigma+d)_k<\alpha(\gamma_1,1-\sigma-d,\sigma+d)_{k+i}$ since we have $(q_i+1)(1-\sigma-d)-(i-q_i-1)(\sigma+d)>0$ using $\frac{1-\{i \sigma\}}{i}>d$. Thus, $\alpha=\alpha(\gamma_1,1-\sigma-d,\sigma+d)$. 

We have proved that if $\alpha$ is a permutation with maximal pattern complexity equal to $n$, and the first string $\gamma_1$ of its table $\gamma$ is a Sturmian word, then $\alpha$ is a Sturmian permutation. It remains to consider the case when $\gamma_1$ is not Sturmian and thus is ultimately periodic.

\section{Proof of uniqueness: periodic case}
We are going to prove that if $\alpha$ is not ultimately periodic and $\gamma_1$ is ultimately periodic, then $p_\alpha^*(n)>n$ for some $n>1$. 

For $n\in \mathbb{N}$, let $\tau^n\alpha$ be the $\mathbb{N}$-permutation such that $(\tau^n\alpha)_i<(\tau^n\alpha)_j$ if and only if $\alpha_{i+n}<\alpha_{j+n}$ for any $i,j\in \mathbb{N}$ with $i\ne j$. Thus, $\tau$ is the {\em shift} on the set of $\mathbb{N}$-permutations. We use the notation $\tau$ also for the shift on the set of words on $\mathbb{N}$. Since the above statement for $\alpha$ follows from that for $\tau^n\alpha$, we'll prove it for $\tau^n\alpha$ such that $\tau^n\gamma_1$ is periodic. Denoting this $\tau^n\alpha$ by $\alpha$, we may assume that $\gamma_1$ is periodic. In the same way, every ultimately periodic sequence defined with respect to $\alpha$ can be consisered as periodic. 

It is convenient to consider arithmetic subpermutations of a permutation $\alpha$. Let us fix a difference $i$ and for each $j=0,\ldots,i-1$ denote by $S_i^j$ the subset $\{ki+j|k \in {\mathbb N}\}$ of $\mathbb N$, called an {\em arithmetic progression} of difference $i$. Now denote by  $\alpha_i^j$ the restriction of $\alpha$ to the set $S_i^j$: $\alpha_i^j=\alpha_{S_i^j}$, and denote by $\alpha_i^{j,k}$ the union of $\alpha_i^j$ and $\alpha_i^k$, that is, the restriction $\alpha_{S_i^j\cup S_i^k}$ of $\alpha$ on $S_i^j\cup S_i^k$. Note that $\alpha$ is not obliged to be an $\mathbb N$-permutation: for all the definitions above, it is sufficient for it to be defined on all values of respective arithmetic progressions.

Let us say that subpermutations $\alpha_i^j$ and $\alpha_i^k$ are {\em adjusted} if $\alpha_i^{j,k}$ is $t'_{j,k}$-periodic for some $t'_{j,k}>0$. (Recall that periodicity was defined for permutations on an arbitrary set, not only for $\mathbb N$-permutations.) Clearly, we can always choose $t'_{j,k}$ divided by $i$, that is, $t'_{j,k}=i t_{j,k}$ for some $t_{j,k}$. It is also clear that to be adjusted with some other subpermutation, a subpermutation must be periodic by itself.

The following lemma has been proved in \cite{ff} in slightly different notation, so we repeat its proof here.
\begin{lemma}\label{l:per}
A permutation defined on a union of infinite arithmetic progressions of difference $i$ is periodic if and only if for all $j,k \in \{0,\ldots,i-1\}$ the subpermutations $\alpha_i^j$ and $\alpha_i^k$ (when well-defined) are adjusted.
\end{lemma} 
{\sc Proof.} The ``only if" part of the proof is obvious since $\alpha_i^{j,k}$ is just restrictions of $\alpha$: if $\alpha$ is $t$-periodic, then so do they.

To prove the ``if" part, we just directly check by the definition that $\alpha$ is $t$-periodic, where $t=i$ lcm $t_{j,k}$, and the lcm (i.e. least common multiple) is taken over all pairs of allowed $j$ and $k$. Indeed, if we take $j'\in S_i^j$ and $k'\in S_i^k$ for any $j$ and $k$, we immediately see that $\gamma_{j'k'}=\gamma_{j'+t'_{j,k},k'+t'_{j,k}}=\gamma_{j'+2t'_{j,k},k'+2t'_{j,k}}=\ldots=\gamma_{j'+t,k'+t}$ which means the $t$-periodicity. \hfill $\Box$

In particular, this lemma holds for all $\mathbb N$-permutations.

Note also that each $i$-periodic permutation consists of $i$ monotonic subpermutations since we have $\gamma_{j,i+j}=\gamma_{i+j,2i+j}=\ldots=\gamma_{ni+j, (n+1)i+j}$ for all $n$.

\begin{claim}\label{perper}
If the maximal pattern complexity of an infinite permutation $\alpha$ satisfies $p^*_{\alpha}(n)=n$, and the sequence $\gamma_1$ is periodic, then for each $i$ and $j$ the sequence $\gamma_i^j$ is periodic.
\end{claim}
{\sc  Proof.} Clearly, if $\gamma_1$ is 1-periodic, then $\alpha$ is monotonic, and there is nothing to be proved. So, we may assume that the minimal period $p$ of $\gamma_1$ is greater than 1, and thus both symbols $<$ and $>$ occur in $\gamma_1$:  moreover, there exist some $k$ and $l$ such that $\gamma_{pn+k,pn+k+1}=<$ and $\gamma_{pn+l,pn+l+1}=>$ for all $n\in \mathbb{N}$.

Due to Lemma \ref{l1}, the sequence $\gamma_i^j$ is either periodic or Sturmian. Suppose it is Sturmian. Then, $\alpha_i^j$ is not periodic, and thus its maximal pattern complexity is at least $n$. The patterns well-defined on $S_i^j$ are exactly those of the form $T=(0,im_1,\ldots,im_n)$ for non-negative $m_1,\ldots,m_n$. Since the maximal pattern complexity of $\alpha_i^j$ cannot be greater than that of $\alpha$, it is equal to $n$. But applying patterns well-defined on $\alpha_i^j$ to $\alpha$ as a whole must not increase the complexity, which immediately means that the language of factors of any subpermutation $\alpha_i^{j'}$ of the same difference $i$ is equal to that of $\alpha_i^{j'}$. In particular, for all $j'$, the sequences $\gamma_i^{j'}$ are Sturmian.

Now consider the pattern $T=(0,1,i)$. By the definition of $k$, for any large $n$ we have that the relation between the first two entries of $\alpha_{T+np+k}$ is 
$<$. At the same time, the relation between $\alpha_{pn+k}$ and $\alpha_{pn+k+i}$ takes both values with different $n$ since positions $k,pi+k,2pi+k,\ldots$ form an arithmetic progression which is a subset of $S_i^k$, and thus elements of the Sturmian word $\gamma_{i}^k$ appearing in this arithmetic progression (of difference $p$ with respect to it) take both values $<$ and $>$. Symmetrically, for any $n$, the relation between the first two entries of 
$\alpha_{T+np+l}$ is 
$>$, 
and the relation between the first and the last elements again takes two values. So, 
$p^*_{\alpha}(3)\geq p_{\alpha}(T)\geq 4$, 
contradicting to our assumption. \hfill $\Box$

\begin{claim}
If the maximal pattern complexity of an infinite permutation $\alpha$ is $p^*_{\alpha}(n)=n$, and the sequence $\gamma_1$ is $p$-periodic, then there exists some $i'$ such that the subpermutation $\alpha_p^{i'}$ is monotonic. 
\end{claim}
{\sc Proof.} 
First of all, we have $p>1$ since otherwise $\alpha$ is monotonic and thus periodic, and its maximal pattern complexity is ultimately constant. Thus, $\gamma_1$ contains both symbols $<$ and $>$: say, the symbols $\gamma_{np+i_1,np+i_1+1}$ for all $n\in \mathbb{N}$ and some $i_1\in \{0,\ldots,p-1\}$ are equal to $<$, and the symbols $\gamma_{np+i_2}$ for all $n\in \mathbb{N}$ and some $i_2\in \{0,\ldots,p-1\}$ are equal to $>$.

Consider the window $T=(0,1,p)$. We must have $p_{\alpha}(T)\leq 3$. Since we always have $\alpha_{T+n_1p+i_1}\neq \alpha_{T+n_2p+i_2}$ for any $n_1,n_2\in \mathbb{N}$, one of the sets $\{\alpha_{T+np+i_1}|n>N\}$ and $\{\alpha_{T+np+i_2}|n>N\}$ (and thus in particular one of the sets $\{\gamma_{np+i_1,n(p+1)+i_1}|n>N\}$ and $\{\gamma_{np+i_2,n(p+1)+i_2}|n>N\}$) is of cardinality one. So, 
either $\alpha_p^{i_1}$ or $\alpha_p^{i_2}$, denoted below by $\alpha_p^{i'}$, is monotonic. \hfill $\Box$

\begin{claim}
If the maximal pattern complexity of an infinite permutation $\alpha$ is $p^*_{\alpha}(n)=n$, and the sequence $\gamma_1$ is periodic, then there exists some $t$ such that all the subpermutations $\alpha_t^{i}$, $i=0,\ldots,t-1$, are monotonic.
\end{claim}
{\sc Proof.} Let $p$ be the minimal period of $\gamma_1$. Consider all the subsequences $\alpha_p^{j}$ with $j=0,\ldots,p-1$. 

Suppose first that some $\alpha_p^{j}$ is $q(j)$-periodic (as a $S_p^{j}$-permutation). Then all its arithmetic subpermutations of difference $q(j)$ are monotonic.


Now consider some of $\alpha_p^{j}$ which is not periodic. However, the word of relations $\gamma_{p}^{j}$ has to be periodic due to the Claim \ref{perper}.
Let us denote its minimal period by $pq$ ($p$ appears here since we consider $\gamma_{p}^{j}$ as a word defined on $S_p^{j}$ not on $\mathbb N$). Since $\alpha_{p}^{j}$ is not monotonic, we have $q\geq 2$, and thus $\gamma_p^{j}$ contains both symbols $<$ and $>$ in the period. 

Clearly, $p_{\alpha_p^{j}}(T_{q+1})\geq p_{\gamma_p^{j}}(T_{q})=q$, where $T_{n}=(0,p,\ldots,(n-1)p)$ for all $n$ and $\alpha_p^{j}$. 

Suppose first that $p_{\alpha_p^{j}}(T_{q+1})>q$. Note that among $T_{q+1}$-factors of $\alpha_p^{j}$, there are no monotonic ones since $\gamma_p^{j}$ is $pq$-periodic and contains both symbols $<$ and $>$ in the period. But $\alpha_{T_{q+1}+i'}$ is monotonic due to the previous claim, and thus 
 $p^*_{\alpha}(q+1)\geq p_{\alpha}(T_{q+1})>q+1$. A contradiction to the minimality of $p^*_{\alpha}$.
 
 Now suppose that $p_{\alpha_p^{j}}(T_{q+1})=q$. This means that each $T_{q+1}$-factor $\alpha_{T_{q+1}+np+j}$ of $\alpha_p^{j}$ for $n\in\mathbb{N}$, and in particular the relation between its first entry $\alpha_{np+j}$ and last entry $\alpha_{(n+q)p+j}$, is determined by the underlying $T_q$-factor 
of $\gamma_p^{j}$ and thus just by the residue of $n$ modulo $q$. So, each of the subsequences $\alpha_{qp}^{n_0p+j}$, where $n_0=0,\ldots,q-1$, is monotonic. Denote $pq=q(j)$. 

Now $q(j)$ is defined for all $j=0,\ldots,p-1$, and all arithmetic subpermutations of $\alpha_p^j$ of difference $q(j)$ are monotonic. Defining $t=$lcm$_{j} q(j)$, we see that all the arithmetic subpermutations of $\alpha$ of difference $t$ are also monotonic, which was to be proved. \hfill $\Box$

So, let $\alpha$ be an infinite permutation such that $p^*_{\alpha}(n)=n$, and the sequence $\gamma_1$ be periodic. Due to the previous Claim and Lemma \ref{l:per}, we see that there exist two 
 subpermutations $\alpha_t^j$ and $\alpha_t^r$ which are monotonic and not adjusted. Without loss of generality we may assume that $j=0$ and both subpermutations are increasing: indeed, if one of them is increasing and the other is decreasing, $\alpha_t^{0,r}$ is $t$-periodic (starting from the point when the subpermutations intersect, if it exists). If they both are decreasing, we just may consider the situation symmetrically.

It is also convenient to denote $\alpha_t^0$ and $\alpha_t^r$ by ($\mathbb N$-permutations) $\chi$ and $\psi$ so that $\chi_i=\alpha_{it}$ and $\psi_i=\alpha_{it+r}$ for all $i \geq 0$. Both permutations are monotonically increasing: $\psi_i<\psi_{i+1}$ and $\chi_i< \chi_{i+1}$ for all $i$.

Note that the fact that $\alpha_t^{0,r}$ is not periodic means in particular that for each $i$ there exists some $v(i)$ such that 
\[\psi_i < \chi_{v(i)}, \mbox{~and~} v(i) \mbox{~is the minimal number with this property.}\]
In particular, if $v(i)>0$, we have $\chi_{v(i)-1}<\psi_i$. 

Symmetrically, for each $j$ there exists some $w(j)$ such that $\chi_j<\psi_{w(i)}$.

Consider first the situation when the modulo $|i-v(i)|$ is bounded: for all $i$, we have $|i-v(i)|<c$.

\begin{lemma}
Permutations $\alpha$ of $p^*_{\alpha}(n)=n$ having periodic sequence $\gamma_1$, not adjusted monotonic subpermutations $\alpha_t^0$ and $\alpha_t^r$, and $|i-v(i)|<c$ for all $i$, do not exist.
\end{lemma}
{\sc Proof.} Suppose such a permutation exists. It follows from the property $|i-v(i)|<c$ for all $i$ that for all $i,n\geq 0$ we have $\psi_i< \chi_{i+c+n}$ and 
$\chi_{i-c-n}<\psi_i$ (of course, the latter inequality is valid only when $i-c-n\geq 0$).

So, for all $n$ we see that all the entries of sequences $\gamma_s$ with $s>(c+n+1)t$ which describe the relations between elements of $\alpha_{t}^{0,r}$ are equal to $<$.

At the same time, we know from Claim \ref{perper} that all sequences $\gamma_s$, and thus their restrictions to $S_{t}^{0}\cup S_t^r$, are periodic. Let us denote the period of $\gamma_s$ by $q_s$; then the restriction of $\gamma_s$ to $S_{t}^{0}\cup S_t^r$ is at most $q_s$-periodic (due to the definition of periodicity involving arbitrary distance between compared periods). Denote by $q$ the least common multiple of all numbers $q_s$ with $s\leq (c+n+1)t$. Then we can check directly that $\alpha_{t}^{0,k}$ is also $q$-periodic, and thus $\alpha_t^0$ and $\alpha_t^r$ are adjusted. A contradiction. \hfill $\Box$

It remains to consider the case when $|i-v(i)|$ is not bounded with $i$: due to the symmetry between $\chi$ and $\psi$, it is sufficient to consider the case when for each $c$ there is some $i$ such that $i-v(i)>c$.

\begin{lemma}
If increasing subpermutations $\alpha_t^0=\chi$ and $\alpha_t^r=\psi$ are not adjusted, and the difference $i-v(i)$ is not bounded with $i$, then $p^*_{\alpha}(4)\geq 5$.
\end{lemma}
{\sc Proof.}
Let us point out a pattern $T$ of length 4 such that $p_{\alpha}(T)\geq 5$. To do it, we need to proof two auxiliary statements.
\begin{claim}
For each $i,j$ there exists some $k$ such that $\psi_{j+k}< \chi_{i+k}$.
\end{claim}
{\sc Proof.} Consider some $n$ such that $n-v(n)>j$, so that $\psi_n<\chi_{v(n)}$ and since both subpermutations are increasing,
$\psi_{j+v(n)}<\psi_n<\chi_{v(n)}\leq \chi_{i+v(n)}$. So, we may take $k=v(n)$. \hfill $\Box$
\begin{claim}\label{cllll}
For each $l$ such that $\chi_1<\psi_l$ there exist some $k_1$ and $k_2$ such that 
\begin{eqnarray*}
\chi_{k_1}<\psi_{l+k_1}<\chi_{k_1+1}\\
\psi_{l+k_2}<\chi_{k_2}<\chi_{k_2+1}.
\end{eqnarray*}
\end{claim}
{\sc Proof.} The number $k_1$ can be found as the minimal number $k$ such that $\psi_{l+k}<\chi_{1+k}$: it exists due to the previous claims, and the fact that it is minimal gives us $\chi_{k_1}=\chi_{1+(k_1-1)}<\psi_{l+k_1-1}<\psi_{l+k_1}$. The number $k_2$ can be found directly from the previous claim as a number such that $\psi_{l+k_2}<\chi_{k_2}$. \hfill $\Box$ 

{\sc Proof of the lemma.} Let us take an arbitrary $l$ such that $\chi_1<\psi_l$, and choose $k_1$ and $k_2$ as described in Claim \ref{cllll}. Now let us choose some $m>l$ such that $\chi_{1+k_1}<\psi_{m+k_1}$ and $\chi_{1+k_2}<\psi_{m+k_2}$ (such $m$ exists since we can take just the greater of the two numbers satisfying these inequations separately).

Let us apply Claim \ref{cllll} to $m$ instead of $l$ and define $k_3$ and $k_4$ so that 
\begin{eqnarray*}
\chi_{k_3}<\psi_{m+k_3}<\chi_{k_3+1},\\
\psi_{m+k_4}<\chi_{k_4}<\chi_{k_4+1}.
\end{eqnarray*}
Also, to unify the notation, suppose that $k_0=0$.
Now consider the 4-window $T=(0,t,lt+r,mt+r)$ and $T$-permutations $\alpha_{T+tk_i}$, where $i=0,1,\ldots,4$. By the definition, for each $i$ the permutation $\alpha_{T+tk_i}$ involves as entries exactly the elements $\chi_{k_i}$, $\chi_{k_i+1}$, $\psi_{l+k_i}$, $\psi_{m+k_i}$. Now it remains to record that all the five permutations $\alpha_{T+tk_i}$ for $i=0,1,\ldots,4$ are different. Indeed, consider the 4-tuples $R_i=(\gamma_{k_it,(k_i+l)t+r}, \gamma_{(k_i+1)t,(k_i+l)t+r}, \gamma_{k_it,(k_i+m)t+r}, \gamma_{(k_i+1)t,(k_i+m)t+r})$ and see that by the construction, $R_0=(<,<,<,<)$, $R_1=(<,>,<,<)$, $R_2=(>,>,<,<)$, 
$R_3=(*,>,<,>)$, and $R_4=(>,>,>,>)$ for some value of $*$. But each $R_i$ contains just a part of information determining $\alpha_{T+tk_i}$. Thus $p^*_{\alpha}(4)\geq p_{\alpha}(T)\geq 5$. \hfill $\Box$

We excluded all possibilities when the maximal pattern complexity of an infinite permutation with the periodic string $\gamma_1$ could be equal to $p^*_{\alpha}(n)\equiv n$. So, the summarizing result of the paper is the following

\begin{theorem}
If an infinite permutation $\alpha$ is not periodic, then $p^*_{\alpha}(n)\geq n$ for any $n$. Moreover, $p^*_{\alpha}(n) \equiv n$ if and only if $\alpha$ is a Sturmian permutation.
\end{theorem}

\end{document}